\documentclass[12pt]{amsart}
\usepackage{amssymb,amsmath,amsthm}

\newtheorem{theorem}{Theorem}[section]
\newtheorem{lem}[theorem]{Lemma}
\newtheorem{prop}[theorem]{Proposition}
\newtheorem{cor}[theorem]{Corollary}
\theoremstyle{definition}
\newtheorem{definition}[theorem]{Definition}
\newtheorem{example}[theorem]{Example}

\newtheorem{convention}[theorem]{Convention}

\newtheorem{alg}[theorem]{Algorithm}

\numberwithin{equation}{section}

\renewcommand{\labelenumi}{(\roman{enumi})}

\def\haa{\hskip 16pt\noindent\hangindent=16pt\hangafter=1}
\def\haaa{\noindent\hskip 64pt\noindent\hangindent=64pt\hangafter=1}

\def\kx     {k[{\mathbf x}]}
\def\kxx        {k({\mathbf x})}
\def\lexp       {{{\mathcal LE}}}
\def\lc         {{{\mathcal LC}}}
\def\lexpz       {{{\mathcal LE}_z}}
\def\lcz         {{{\mathcal LC}_z}}
\def\Qset   {\mathbb Q}
\def\Zset   {\mathbb Z}

\def\N   {\mathbb N}

\def\kq         {k\langle \langle t^\Qset \rangle\rangle}
\def\kcq         {\overline{k}\langle \langle t^\Qset \rangle\rangle}
\def\Q      {\mathbb{Q}}
\def\Z      {\mathbb{Z}}
\def\supp   {\mbox{Supp}}

\def\gcd    {\mbox{gcd}}
\def\lcm    {\mbox{lcm}}

\def\kx         {{k[\mathbf{x}]}}
\def\kxx        {{k(\mathbf{x})}}
\def\lexp       {{{\mathcal LE}}}

\def\lcm    {\mbox{\rm lcm}}

\def\kq         {k\langle\langle t^\Q\rangle\rangle}
\def\skl        {k(( t^{-1}))}
\def\supp       {\mbox{\rm Supp}}
\def\gal        {\mbox{\rm Gal}}

\begin{document}
\bibliographystyle{plain}
\title[Value Monoids]
{Value Monoids of Zero-Dimensional Valuations of Rank One}
\author{Edward Mosteig}
\address{Department of Mathematics, Loyola Marymount University,
Los Angeles, California 90045} \email{emosteig@lmu.edu}

\begin{abstract}
Classically, Gr\"obner bases are computed by first prescribing a
set monomial order.  Moss Sweedler suggested an alternative and developed a framework to perform such computations by
using valuation rings in place of monomial orders.  We build on these
ideas by providing a class of valuations on $k(x,y)$ that
are suitable for this framework.
For these valuations, we compute $\nu(k[x,y]^*)$ and use this to perform computations
concerning ideals in the polynomial ring $k[x,y]$.
Interestingly, for these valuations, some ideals have
a finite Gr\"obner
basis with respect to the valuation that is not a Gr\"obner basis with respect to any
monomial order,
whereas other ideals only have Gr\"obner bases that are infinite with respect to the valuation.
\end{abstract}

\maketitle


\section{Introduction}\label{introduction}

Unless stated otherwise, $k$ will denote an arbitrary field, and
 $\N$ will denote the set of nonnegative integers.
Whenever $R$ is a ring or monoid, we denote by
 $R^*$
the nonzero elements of $R$.

One of the
fundamental ideas of the theory of Gr\"obner bases is that monomial orders are well-orderings on the set
of monomials, which leads us to a natural reduction process using multivariate polynomial division.
In this section, we provide a brief account of a generalized
theory of  Gr\"obner
bases that uses valuations in place of monomial orders, which will
yield a more general reduction process.  The development of this theory can be found in the unpublished
manuscript \cite{sweedler} of Sweedler, and it is briefly discussed in this section solely for the sake of completeness.  In that manuscript, Sweedler
develops the theory in terms of valuation rings. Here we present the same results  in terms of
valuations rather than valuation rings.  Proofs are
omitted since they can all  be found in \cite{sweedler}.

Suppose $k$ is a subfield of a field $F$.
A {\bf valuation on} $F$ is a homomorphism $\nu$ from the additive
group of nonzero elements of $F$ to an ordered group (called the {\bf value group}) such that for $f,g \in F^*$ where $f+g
\not=0$, $\nu(f+g) \le \max\{\nu(f), \nu(g)\}$.  Note that the
triangle inequality was chosen to be opposite of the most common
definition, which is so that our results most closely coincide
with those concerning monomial orders.   For more
details, see \cite{ms1}, \cite{ms2}, and \cite{ms3}.
A {\bf valuation on $F$ over $k$} is a valuation on $F$ such that its restriction
to $k^*$ is the zero map.
For our purposes, we restrict our attention to valuations on rational function fields.  In this setting, we require
that our valuations have the additional properties given in the following definition.

\begin{definition} \label{def:suitableval}
We say that a valuation $\nu$ on $\kxx$ over $k$ is {\bf suitable  relative to } $\kx$
if satisfies the following three properties.
\begin{enumerate}
\item[(i)] For all $f \in k[\mathbf{x}]$, $\nu(f) = 0$ iff $f\in k$.
\item[(ii)] If $\nu(f) = \nu(g)$ where $f,g \in k({\mathbf x})^*$, then
$\exists! \lambda \in k^*$ such that  $f = \lambda g$ or $\nu(f- \lambda g) < \nu(f)$.
\item[(iii)] $\nu(k[\mathbf{x}]^*)$ is a well-ordered monoid.
\end{enumerate}
\end{definition}

When using monomial orders, one must determine divisibility among monomials.
 The analogue for valuations uses arithmetic in the monoid $\nu(\kx^*)$.

\begin{definition}
\label{def:division}  Let
$\nu$ be a valuation on $k({\mathbf x})$. Given $f,
g \in \kx$, we say that $\nu(g)$ {\bf divides} $\nu(f)$, denoted
$\nu(g)\ |\ \nu(f)$, if there exists $h \in
\kx$ such that $\nu(f) = \nu(gh)$.  We say that $h$ is an {\bf
approximate quotient} of $f$ by $g$ (relative to $\nu$), if
$f=gh$, or if $f\not=gh$ and $\nu(f-gh) < \nu(f)$.
\end{definition}

The following simple proposition follows  from the definition above.

\begin{prop}
\label{prop:approxquotient}
 Let $\nu$ be a valuation on $\kxx$ over $k$ that
is suitable relative to $\kx$.
Let $f, g \in \kx$.  Then $\nu(g)$
divides $\nu(f)$ if and only if there exists  an approximate
quotient $h$  of $f$ by $g$.
\end{prop}

The following is a generalized form of the standard polynomial reduction
algorithm that makes use of  valuations.

\begin{alg}
\label{alg:genred}
 Let $\nu$ be a valuation on $\kxx$ over $k$ that
is suitable relative to $\kx$.  Let $\nu$ be a valuation on $\kxx$ over $k$.
 Let $I$ be an ideal in
$\kx$ and $G$ be a generating set for $I$.
  The following
algorithm computes a reduction of a polynomial $a\in \kx$ over $G$ relative to $\nu$.

\haa $\bullet$  Set $i=0$ and $f_0 = f$.

\haa $\bullet$ While $f_i \not=0$ and $\nu(g) \ | \ \nu(f_i)$ for some $g \in G$ do:

\haaa  Choose $g_i \in G$ such that $\nu(g_i) \ | \ \nu(f_i)$.
 Let  $h_i$ be an  an approximate quotient of $f_i$ by $g_i$.
 Set $f_{i+1} = f_i - g_ih_i$.  Increment $i$ by 1.

\end{alg}

We say that $f_n$ is the $n$th reductum of $f$ over $G$.  We
say that $f$ reduces to $b$ if $b$ is a reductum of $f$.
It can be shown that if $\nu$ is suitable with respect to $\kx$, then reduction of any element of $\kx$ over $G$ terminates after a finite number of steps.
We will call a subset $G \subset I^*$ a {\bf Gr\"obner basis for $I$ with respect to $\nu$} if it satisfies the equivalent conditions of the
following proposition.

\begin{prop}
\label{prop:gbconds} Let $\nu$ be a valuation on $\kxx$ over $k$ that
is suitable relative to $\kx$.
  Let $I$ be an
ideal in $\kx$ and $G \subseteq I^*$.  The following are equivalent:
\renewcommand{\labelenumi}{(\roman{enumi})}\begin{enumerate}\renewcommand{\labelenumi}{(\roman{enumi})}
\item Every nonzero
element of $I$ has a first reductum over $G$.
\item Every element of $I$ reduces to $0$ over $G$.
\item  Given $f \in \kx$, $f \in I$ if and only if $f$ reduces to $0$ over $G$.
\end{enumerate}
\end{prop}

We can use Gr\"obner bases in the generalized setting to solve the
ideal membership problem in much the same way that we do in the
case of monomial orders.
Just as in the classical case, it can be shown that a Gr\"obner basis with respect to a valuation necessarily
generates the given ideal.  To compute Gr\"obner bases, we must work with ideals of $\nu(\kx^*)$, where
an ideal $J$ of a commutative monoid $M$ is a subset $J \subset M$ such that
for any $m\in M, j \in j$, $j+m \in J$.  The smallest ideal containing $m_1, \dots, m_\ell$
will be denoted $\langle m_1, \dots, m_\ell \rangle$ and is called
the ideal generated by $m_1, \dots, m_\ell$

\begin{definition}
 Let $\nu$ be a valuation on $\kxx$ over $k$ that
is suitable relative to $\kx$.  We say that  $T \subseteq
\nu(\kx^*)$ is {\bf an ideal generating set for $f$
and $g$ with respect to $\nu$} if $T$
generates the ideal $\langle \nu(f) \rangle \cap \langle \nu(g)
\rangle$ in $\nu(\kx^*)$. It can be shown that for
 each $t\in T$ there are $a, b
\in \kx^*$ such that $\nu(af) = \nu(bg) = t$ and
$af=bg$ or $af\not=bg$ and $\nu(af-bg) < t$. This
gives a map $T \rightarrow \kx$, $t \mapsto af-bg$.  The image
of this map is a {\bf syzygy family for $f$ and $g$ indexed by
$T$}. We say that $af-bg$ is the element of the family
corresponding to $t$.
\end{definition}

This definition shows one of the main
differences between the generalized theory using valuations and the
classical theory using monomial orders, namely, that each pair of polynomials may
have many minimal syzygies.  Sweedler constructs an
example in \cite{sweedler} where this family consists of multiple elements.
Using syzygy families, the algorithm below provides a method for constructing a Gr\"obner basis for an nonzero ideal $I$ with
generating set $G$.

\begin{alg} [Gr\"obner Basis Construction Algorithm]
\label{alg:gengbconstruct}  Let $\nu$ be a valuation on $\kxx$ over $k$ that
is suitable relative to $\kx$, and
$G \subseteq I^*$ is a generating set for a nonzero ideal $I$.
\renewcommand{\labelenumi}{(\roman{enumi})}\begin{enumerate}\renewcommand{\labelenumi}{(\roman{enumi})}
\item Set $G_0 = G$.
\item For each pair of distinct elements $g, h \in G$, find a monoid generating set $T^0_{g,h}$
for $g,h$ and a syzygy family $S^0_{g,h}$ for $g,h$ indexed by
$T^0_{g,h}$.  Define $U = \bigcup_{g\not=h \in G} S^0_{g,h}$.
\item  Determine the set $H_i$ of nonzero final reductums that occur from reducing the lements of $U_i$ over $G_i$.
\item If $H_i$ is empty, stop.
\item Define $G_{i+1} = G_i \cup H_i$.
\item For each pair of distinct element $g \in G_{i+1}, h \in H_i$,
find a monoid generating set $T^{i+1}_{g,h}$ for $g,h$ and a syzygy family $S^{i+1}_{g,h}$ for $g,h$ indexed by $T^{i+1}_{g,h}$.  Define $U = \bigcup_{g\not=h \in G} S^{i+1}_{g,h}$.
\item Increment $i$ by 1 and go to step (iii).
\end{enumerate}
\end{alg}

Sweedler shows that if $G$ is finite and
$\nu(I^*)$ is Noetherian (i.e., every ascending chain of ideals
stabilizes), then the construction algorithm can be completed so that it terminates
with a finite Gr\"obner basis.  However, even if $\nu(I^*)$ isn't Noetherian,
the set $\cup_{n=1}^\infty G_n$ is a Gr\"obner basis.

These algorithms will allow us to
 compute Gr\"obner bases using a class of valuations
on $k(x,y)$ originally studied by Zariski in \cite{zariski}.
In Section \ref{valmonoid}, we develop the background necessary to work with
a valuation $\nu$ of this type, and we state one of the main results of the paper, which is
an explicit formula for $\nu(k[x,y]^*)$.
In Section \ref{associatedsequences}, we prove some intermediate results concerning
sequences associated with the valuations developed in Section \ref{valmonoid}.  In
particular, recursive formulas are given for a generating set of $\nu(k[x,y]^*)$.
In Section \ref{repsmonoid}, we build on these ideas to show that certain
 elements of
$\nu(k[x,y]^*)$ have unique representations, which leads to a complete
description of $\nu(k[x,y]^*)$ in Section \ref{construction}.  Finally,
in Section \ref{algorithms}, we
use this description to make the algorithms developed by Sweedler constructive.
With the exception of Section \ref{repsmonoid}, all of the proofs herein are fairly elementary.

\section{Value Groups and Monoids from Power Series}\label{valmonoid}

In this section, we examine a class of valuations of $k(x,y)$
studied by Zariski in \cite{zariski}.   The value groups of these
valuations were explicitly constructed by MacLane and Schilling in
\cite{mac-schi}.  In this section, we state one of our main
results, which is an explicit construction of the restriction of
such valuations to the underlying polynomial ring $k[x,y]$.  Since
the valuations of interest are constructed using generalized power
series, we begin with a review of the relevant concepts.

 We say that a set  $T \subset \Q$ is {\bf
Noetherian} if every subset of $T$ has a largest element.   Given
a function $z: \Q \to k$, the {\bf support} of $z$ is defined by
 $\supp(z) = \{ q \in \Q \mid z(q) \not= 0 \}.$
The collection of {\bf Noetherian power series}, denoted by $\kq$,
consists of all functions from $\Q$ to $k$ with Noetherian
support.  More commonly in the literature, generalized power
series are defined as functions with well-ordered support, and we
will freely use the analogues of these results for Noetherian
power series.  We choose the supports of our series to be opposite
of the usual definition so that our results more closely fit with
the theory of monomial orders and Gr\"obner bases.

As demonstrated in \cite{hahn}, the collection of Noetherian power
series forms a field in which addition is defined pointwise and
multiplication is defined via convolution; i.e., if $z_1, z_2 \in
\kq$ and $q \in \Q$, then $(z_1+z_2) (q) = z_1(q) + z_2(q)$ and
$(z_1z_2)(q)  = \sum_{ u+v=q} z_1(u)  z_2(v)$.  We often write
power series as formal sums: $z = \sum_{s\in{{\mbox{\scriptsize\rm
Supp}}}(z)} z(s)t^s$, where $z(s)$ denotes the image of $s$ under
$z$.

\begin{example}\label{ex:sumprod}
Given the series $z_1  =  t^{1/2} + t^{1/4} + t^{1/8} + \cdots$
and $z_2  =  3t+ 1,$ their sum and product are  $$z_1+ z_2 = 3t +
( t^{1/2} + t^{1/4} + t^{1/8} + \cdots) + 1$$ and $$z_1z_2 =
(3t^{3/2} + 3t^{5/4} + 3t^{9/8} + \cdots ) + (t^{1/2} + t^{1/4} +
t^{1/8} + \cdots).$$
\end{example}

Given a series $z \in \kq$, define the {\bf leading exponent} of
$z$ to be the rational number given by $\lexp (z) = \max\{ s \mid s
\in \supp(z) \}$.     If $s= \lexp(z)$, we denote $z(s)$ by
$\lc(z)$ and call it the {\bf leading coefficient } of $z$. Note
that $\lexp(z_1  z_2) = \lexp(z_1) + \lexp(z_2)$ and $\lc(z_1z_2)
= \lc(z_1)\lc(z_2)$. Moreover, we have $\lexp(z_1 + z_2) \le
\max(\lexp(z_1), \lexp(z_2))$, with equality holding in case
$\lexp(z_1) \not= \lexp(z_2)$.

We say that a nonzero series  $z \in \kq$ is {\bf simple} if it
can be written in the form \begin{equation*}\label{def:simple} z =
\sum_{i=1}^{n} c_i t^{e_i},
\end{equation*}
where $c_i \in k^*, n \in \N^* \cup \{ \infty \}, e_i \in \Q, e_i
>  e_{i+1}$.   Whenever we write a series in this form, we
implicitly assume that each $c_i$ is nonzero and the exponents are
written in descending order. We call ${\mathbf e} = (e_1, e_2,
\dots)$ the {\bf exponent sequence} of $z$.  Now write $e_i =
n_i/d_i$ where $d_i>0$ and $\gcd(n_i,d_i)=1$.
 We define $r_0 =1$ and for $i \ge 1$,
 set $r_i = \lcm(d_1, \dots, d_{i})$ and call  ${\mathbf{r}} = (r_0, r_1, r_2,  \dots)$ the
{\bf ramification sequence} of $z$.

\begin{example}
Consider the simple series $$z=2t^{1/2}+3t^{1/3}+4t^{1/4} +
5t^{1/5} + \cdots.$$ Here $\lexp(z) = 1/2$ and $\lc(z) = 2$.  The
series $z$ has exponent sequence $(1/2, 1/3, 1/4, 1/5, \cdots)$
and ramification sequence $(1,2,6,12,60, \cdots)$.
\end{example}

We are now in a position
to define valuations on $k(x,y)$ based on Noetherian power series.
Let $z\in \kq$ be a Noetherian power series such that $t$ and $z$
are algebraically independent over $k$. Consider the embedding
 $\varphi_z : k(x,y) \to
\kq$, $x \mapsto t$, $y \mapsto z$.  It can be shown that $\lexp$
is a valuation on $\kq$, and hence the composite map $\lexp \circ
\varphi_z : k(x,y) \to \Q$ is a valuation on $k(x,y)$.
Given a valuation $\nu$ on $\kxx$,  $V = \{ f\in \kxx^* \mid \ \nu(f) \le 0 \}$ is
a valuation ring with maximal ideal ${\mathfrak m} = \{ f\in \kxx^* \mid \ \nu(f) < 0 \}$,
in which case $\dim_k (V/{\mathfrak m})$ is the {\bf dimension} of the valuation.
The {\bf rank} of the valuation $\nu$ is defined to be the number of isolated
subgroups of $\nu(\kxx^*)$.  It follows that $\lexp \circ \varphi_z$ is
a zero-dimensional valuation of rank one.

\begin{example}
Let $k$ be a field such that char $k \not=2$. Given $z = t^{1/2} +
t^{1/4} + t^{1/8} + \cdots$,
\begin{eqnarray*}
(\lexp \circ \varphi_z) (x) & = & \lexp(t) = 1 \\
(\lexp \circ \varphi_z) (y) & =&  \lexp(z) = 1/2 \\
(\lexp \circ \varphi_z) (y^2-x) & = & \lexp(z^2-t) = \lexp(
((t+2t^{3/4} +  2t^{5/8} + \cdots ) - t) = 3/4
\end{eqnarray*}
\end{example}

MacLane and Schilling proved the following result in
\cite{mac-schi}:

\begin{theorem}
Let $z\in \kq$ be a simple series such that $t$ and $z$
are algebraically independent over $k$.   If $\mathbf e$ is the
exponent sequence of $z$, then the value group of $\lexp \circ
\varphi_z$ is
$$(\lexp \circ \varphi_z) (k(x,y)^*) = \Z + \Z e_1 + \Z e_2 + \cdots $$
\end{theorem}

One of the primary goals of this paper is to restrict the
valuation to the polynomial ring $k[x,y]$ and compute
\begin{equation}\label{eq:mz}
{\Lambda_{}} = (\lexp \circ \varphi_z) (k[x,y]^*) = \{ \lexp(f(t,z)) \mid
f(x,y) \in k(x,y)^* \},
\end{equation}
which we call the {\bf value monoid with respect to $z$}.

Now suppose $z$ is a simple series with exponent sequence
${\mathbf e}$ and ramification sequence ${\mathbf r}$. The
sequence obtained from the ramification sequence $\{ r_i \}_{i \in
\N}$  by removing repetitions is called the {\bf reduced
ramification sequence} and is denoted $\{ r_i^{red} \}_{i \in
\N}$.  For each $i \in \N$, denote by $l(i)$ the smallest natural
number such that $r_i^{red} = r_{l(i)}$; i.e.,
\begin{equation}\label{eq:ls}
l(i) = \min \{ j \in \N \mid r_j = r_i^{red} \}.
\end{equation}

\begin{example}
The series
 $$z= t^2 + t^{3/2} + t^{1/2} + t^{1/3} + t^{1/5} + t^{1/7} +
t^{1/11} + \cdots$$ has ramification sequence
$${\mathbf{r}} = (1,2,2,6,30,210,2310, \dots),$$
and hence has reduced ramification sequence
$$ (1,2,6,30,210,2310, \dots).$$
Thus $l(0) = 0, \ l(1) = 1, \  l(i) = i+1  \mbox{ for } i \ge 2.$
\end{example}

 We define the
{\bf bounding sequence} ${\mathbf u} = (u_0, u_1, u_2, \dots)$
given by $u_0 = 0$, and for $i \ge 1$,
\begin{equation}\label{eq:u}
u_i = \sum_{j=0}^{i-1} \Big( \frac{r_i}{r_{j}} -
\frac{r_i}{r_{j+1}} \Big) e_{j+1}.
\end{equation}
 For $i \ge 1$, we define the {\bf monoid generating sequence}:
\begin{eqnarray}\label{def:rho}
\rho_i & = & u_{l(i)-1} + e_{l(i)}.
\end{eqnarray}
We can fully describe the value monoid with respect to $z$ in
terms of the monoid generating sequence.
  The following result will
be proved in Section \ref{construction}
(in fact, it follows directly from the stronger result given in
Theorem \ref{Mzunique}).

\begin{theorem}\label{thm:valmonoid}
Let $z\in \kq$ be a simple series such that $t$
and $z$ are algebraically independent over $k$.  Assume further
that the components of the exponent sequence are positive and no
component is divisible by the characteristic of $k$.
 Then the value monoid with respect to $z$ is
$$ \Lambda = (\lexp \circ \varphi_z) (k[x,y]^*) = \N + \N \rho_1 + \N \rho_2 + \cdots $$
\end{theorem}

It is of interest to determine whether this result can be
generalized.  In particular, it would be nice to compute the value monoid
after either removing the restriction that the exponent sequence must be
positive or permitting some of the
components of the exponent sequence to be divisible by the
characteristic of the the ground field.

\section{Associated Sequences}
\label{associatedsequences}

In this section, we prove some elementary results about the sequences described in the previous
section.  In particular, we will construct recurrence relations and formulas concerning the
monoid generating sequence.  To this end,
there is one more sequence that will be needed in the sequel.
Using the ramification sequence $\mathbf r$ of a simple series $z$
 and the formula (\ref{eq:ls}), we define
 {\bf partial ramification sequence} by
\begin{eqnarray*}
s_i &  = &  r_{l(i)}/r_{l(i-1)} = r_{l(i)}/
r_{l(i)-1}\label{eq:s}. \end{eqnarray*}

\begin{center}
\framebox{
\parbox{6in}
{
\begin{convention}\label{conventions}
For the remainder of this paper, we adopt the following
conventions.
\begin{itemize}
\item The series $z \in\kq$ is simple with positive support.
\item The series $z$ is transcendental over $k(t)$.
\item The value monoid of $z$ is denoted ${\Lambda_{}}$. \item The exponent sequence of $z$ is denoted
${\mathbf e} = (e_1, e_2, e_3, \dots)$.
\item No component of the exponent sequence is divisible by char $k$.
\item The ramification
sequence of $z$ is denoted ${\mathbf r} = (r_0,r_1,r_2, \dots)$.
\item The bounding sequence of $z$ is denoted ${\mathbf u} = (u_0,
u_1, u_2, \dots)$.
\item The function $l(i)$ is defined in  (\ref{eq:ls}).
\item The monoid generating
sequence of $z$ is denoted ${\mathbf \rho} = (\rho_1, \rho_2,
\rho_3, \dots)$.
 \item The partial
ramification sequence of $z$ is given by ${\mathbf s} = (s_1, s_2,
s_3, \dots)$.
\end{itemize}
\end{convention}
}}
\end{center}

Since $l(i)$ marks the index where the ramification index
increases, we have $r_j = r_{l(i)}$ for  $l(i) \le j < l(i+1)$,
and so
\begin{equation}\label{eq:rj}
r_{j}/r_{j-1} = 1 \mbox{ \ for \ } l(i) < j < l(i+1).
\end{equation}
In particular, this yields
\begin{equation}\label{eq:redram}
r_{l(i-1)} = r_{l(i) - 1}
\end{equation}
 and
\begin{equation}\label{eq:redu}
u_{l(i-1)} = u_{l(i) - 1}
\end{equation}
despite the fact that $e_{l(i-1)}$ and $e_{l(i)-1}$ need not be
the same.

Note that the ramification sequence of a series $z\in\kq$
increases without bound unless $z\in k((t^{1/n}))$ for some
$n\in\N$.  However, it is still possible that the ramification
sequence occasionally (even infinitely many times) stabilizes for
a finite number of steps. Whenever the ramification sequence
stabilizes for a number of indices, the sequence
$\{u_i\}_{i\in\N}$ also stabilizes, as seen in the next result.

\begin{lem}\label{lem:staticuandr}
If $r_i=r_k$ for indices $i$ and $k$, then $u_i=u_k$.
\end{lem}

\begin{proof}
The result is trivial if $i=k$, so we assume $i<k$. Since
$r_i=r_k$, it follows that $r_j = r_{j+1}$ for $i \le j \le k-1$,
and so by (\ref{eq:u}),
$u_k =  \sum_{j=0}^{k-1} \Big( \frac{r_k}{r_{j}} -
\frac{r_k}{r_{j+1}} \Big) e_{j+1} = \sum_{j=0}^{k-1} \Big( \frac{r_i}{r_{j}} -
\frac{r_i}{r_{j+1}} \Big) e_{j+1} = u_i + \sum_{j=i}^{k-1} \Big(
\frac{r_i}{r_{j}} - \frac{r_i}{r_{j+1}} \Big) e_{j+1}
 = u_i.$
\end{proof}

Since our main objective is to prove that ${\Lambda_{}}$ is generated by
the sequence $1, \rho_1, \rho_2, \dots$, we must first justify
some elementary properties that allow us to understand better the behavior of this sequence.
We begin by showing that the monoid generating sequence satisfies a simple recursive relation.

\begin{lem}\label{lem:therecurrence}
The monoid generating sequence in (\ref{def:rho}) satisfies the following recurrence
relation:
\begin{eqnarray*}
\rho_1 & = & e_{l(1)}; \\
\rho_{i+1} & = & s_i \rho_i - e_{l(i)} + e_{l(i+1)}.
\end{eqnarray*}
\end{lem}

\begin{proof}
By (\ref{eq:redu}), $u_{l(1)-1} = u_{l(1-1)} = u_0 = 0$, and so by
(\ref{def:rho}), $\rho_1 = e_{l(1)}$. Also, we have by
(\ref{eq:u}),
\begin{eqnarray*}
 u_{m} + e_{m+1}
 & = & \sum_{j=0}^{m-1} \left(
\frac{r_{m}}{r_j} - \frac{r_{m}}{r_{j+1}}
\right) e_{j+1} + e_{m+1}\\
& = & \left(\frac{r_{m}}{r_{m-1}}\right) \sum_{j=0}^{m-2} \left(
\frac{r_{m-1}}{r_j} - \frac{r_{m-1}}{r_{j+1}} \right) e_{j+1}
+\left( \frac{r_{m}}{r_{m-1}} -\frac{r_m}{r_m} \right) e_{m}
+  e_{m+1}\\
& = & \left(\frac{r_{m}}{r_{m-1}}\right) \sum_{j=0}^{m-2} \left(
\frac{r_{m-1}}{r_j} - \frac{r_{m-1}}{r_{j+1}} \right) e_{j+1}
+\left( \frac{r_{m}}{r_{m-1}} \right)e_{m}
-  e_{m} +  e_{m+1}\\
& = & \left(\frac{r_{m}}{r_{m-1}}\right)  [u_{m-1} + e_{m}] -e_{m}
+ e_{m+1},
\end{eqnarray*}
and so \begin{equation}\label{eq:gammaprelim} \gamma_{m+1} =
\left( \frac{r_m}{r_{m-1}} \right) \gamma_m - e_m + e_{m+1}
\end{equation}
where $\gamma_m := u_{m-1} + e_m$.  Replacing $m$ by $l(i)$, we
obtain
\begin{equation}\label{eq:gammali} \gamma_{l(i)+1}
= \left( \frac{r_{l(i)}}{r_{l(i)-1}} \right) \gamma_{l(i)} -
e_{l(i)} + e_{l(i)+1} = s_i \gamma_{l(i)} - e_{l(i)} + e_{l(i)+1}
.\end{equation} If $l(i) < m < l(i+1)$, then $r_m/r_{m-1} = 1$ by
(\ref{eq:rj}), and so (\ref{eq:gammaprelim}) yields
$$\gamma_{m+1} = \gamma_m - e_m + e_{m+1}.$$
Multiple applications of this formula yields a telescoping sum, and
so
\begin{eqnarray*}
\gamma_{l(i+1)} &  = &  \gamma_{l(i+1)-1} - e_{l(i+1)-1} +
e_{l(i+1)} \\
 &  = &  (\gamma_{l(i+1)-2}- e_{l(i+1)-2}  + e_{l(i+1) -1} ) - e_{l(i+1)-1} +
e_{l(i+1)} \\
& = &  \gamma_{l(i+1)-2}- e_{l(i+1)-2}  +
e_{l(i+1)} \\
&  &  \vdots\\
& = &  \gamma_{l(i)+1}- e_{l(i)+1}  + e_{l(i+1)}.
\end{eqnarray*}
This equation in conjunction with (\ref{eq:gammali}) yields
\begin{eqnarray*}\label{eq:gammar}
\gamma_{l(i+1)} & =  & \gamma_{l(i)+1}- e_{l(i)+1}  + e_{l(i+1)}\\
& = &  s_i \gamma_{l(i)} - e_{l(i)} + e_{l(i)+1} -
e_{l(i)+1}  + e_{l(i+1)} \\
& = &  s_i \gamma_{l(i)} - e_{l(i)}   + e_{l(i+1)},
\end{eqnarray*}
and since $\rho_i = u_{l(i)-1} + e_{l(i)} = \gamma_{l(i)}$ for all
$i$, we have
$$\rho_{i+1} = s_i \rho_i - e_{l(i)} + e_{l(i+1)}.$$
\end{proof}

We can also construct a recursive formula for the terms of the ramification sequence, as given in the next
result.

\begin{lem}\label{lem:sumram} For $i \in \N$,
\begin{equation*}
r_{l(i)}  =  1 + \sum_{j=1}^i (s_j-1) r_{l(j-1)}.
\end{equation*}
\end{lem}

\begin{proof}
This follows from the simple computation $$\sum_{j=1}^{i}
(s_j-1)r_{l(j-1)} = \sum_{j=1}^{i}
((r_{l(j)}/r_{l(j-1)})-1)r_{l(j-1)} = \sum_{j=1}^{i}
(r_{l(j)}-r_{l(j-1)}) = r_{l(i)} - r_{l(0)}= r_{l(i)} - 1.
$$
For the case $i=0$, we take the summation $\sum_{j=1}^i (s_j-1) r_{l(j-1)}$ to be
$0$.
\end{proof}

Using Lemma \ref{lem:therecurrence}, we can construct yet another recurrence relation for
the terms of the monoid generating sequence.

\begin{lem}\label{lem:difflams}
For $i \ge 1$, \begin{equation*}\label{eq:rhodiff}\rho_i =
\sum_{j=i}^{i-1}(s_j-1) \rho_j + e_{l(i)}.\end{equation*}
\end{lem}

\begin{proof}
We proceed by induction.  If $i=1$, then by Lemma
\ref{lem:therecurrence},
$$\rho_1  = e_{l(1)} = 0 + e_{l(1)} = \sum_{j=1}^0
(s_j-1)\rho_j + e_{l(1)}$$ since the summation that appears is
empty. Now suppose the statement holds for the index $i$. Then by
Lemma \ref{lem:therecurrence} and the induction hypothesis,
\begin{eqnarray*}
\rho_{i+1} - \sum_{j=1}^{i} (s_j - 1) \rho_j  & = & \rho_{i+1} -
(s_i-1) \rho_i - \sum_{j=1}^{i-1}
(s_j-1)\rho_j \\
& = & \rho_{i+1} - (s_i-1) \rho_i - (\rho_i - e_{l(i)}) \\
& = & s_i\rho_i - e_{l(i)} + e_{l(i+1)} - (s_i-1) \rho_i
- (\rho_i - e_{l(i)}) \\
& = & e_{l(i+1)}.
\end{eqnarray*}
\end{proof}

Using this lemma, we can extract information about the denominators of the
components of the monoid generating sequence, as shown in the next three results.
Given $q\in \Q$, $q\Z$ denotes the set $\{qz \mid z \in \Z\}$.

\begin{lem}\label{cor:rhoresidue}
For $i \ge 1$, $\rho_i \in (1/r_{l(i)}) \Z - (1/r_{l(i-1)}) \Z$.
\end{lem}

\begin{proof}
The result follows by a simple induction.  Indeed, $\rho_1 =
e_{l(1)} \in  (1/r_{l(1)}) \Z -
\Z = (1/r_{l(1)}) \Z  - (1/r_{l(0)})\Z $.  Now, assuming that $\rho_i \in (1/r_{l(i)}) \Z$, we see by
Lemma \ref{lem:difflams}, $\rho_{i+1} = \sum_{j=1}^{i}(s_j-1)
\rho_j + e_{l(i+1)}$.  Since $\rho_j \in (1/r_{l(j)}) \Z \subset
(1/r_{l(i)}) \Z$ for $1 \le j \le i$, we have
$\sum_{j=1}^{i}(s_j-1) \rho_j \in (1/r_{l(i)}) \Z$. Moreover,
$e_{l(i+1)} \in (1/r_{l(i+1)}) \Z -(1/r_{l(i)}) \Z $, and so
$\rho_{i+1}
 \in (1/r_{l(i+1)})
\Z - (1/r_{l(i)})
\Z .$
\end{proof}

\begin{lem}\label{lem:rhoalpha}
If we write $\rho_i = {c_i}/r_{l(i)}$, then $\gcd({c_i},s_i)=1$.
\end{lem}

\begin{proof}
Rewrite the expression $\rho_i = {c_i} /r_{l(i)}$ in lowest terms:
$\rho_i = {\alpha_i}/{\beta_i}$, ${\alpha_i},{\beta_i} \in \N^*$
where $\gcd({\alpha_i},{\beta_i}) =1$. Then ${c_i} = {\alpha_i}
r_{l(i)}/{\beta_i} = {\alpha_i}
\lcm(r_{l(i-1)},{\beta_i})/{\beta_i} = {\alpha_i}
r_{l(i-1)}/\gcd(r_{l(i-1)},{\beta_i}).$ Also $r_{l(i)}/r_{l(i-1)}
= \lcm(r_{l(i-1)},{\beta_i})/r_{l(i-1)} =
{\beta_i}/\gcd(r_{l(i-1)},{\beta_i}).$ Therefore,
$$\gcd({c_i},r_{l(i)}/r_{l(i-1)}) = \gcd({\alpha_i}
r_{l(i-1)}/\gcd(r_{l(i-1)},{\beta_i}),{\beta_i}/\gcd(r_{l(i-1)},{\beta_i})).$$
Since $\gcd({\alpha_i,\beta_i}) = 1$ and
$\gcd(r_{l(i-1)}/\gcd(r_{l(i-1)},{\beta_i}),
{\beta_i}/\gcd(r_{l(i-1)},{\beta_i}))=1$, we have $\gcd(c_i,s_i) =
\gcd({c_i},r_{l(i)}/r_{l(i-1)})=1$.
\end{proof}

\begin{lem}\label{lem:lambdaresidue}
If $0 \le d_j < s_j$ for $1 \le j \le i$ and $d_i \not=0$, then
\begin{equation}\label{lambdaresexpress}
\sum_{j=1}^i d_j \rho_j \in (1/r_{l(i)})\Z - (1/r_{l(i-1)})\Z.
\end{equation}
\end{lem}

\begin{proof}
For $j \le i$, we have by Lemma \ref{cor:rhoresidue}, $\rho_j
\in (1/r_{l(j)}) \Z \subset (1/r_{l(i)}) \Z$, and so $\sum_{j=1}^i
d_j \rho_j \in (1/r_{l(i)})\Z$.  We now must prove $\sum_{j=1}^i
d_j \rho_j \not\in (1/r_{l(i-1)})\Z$ by induction.

First, we show that $d_j \rho_j \not\in (1/r_{l(j-1)}) \Z $
whenever $0 < d_j < s_j$.
 Write $\rho_j = c_j/r_{l(j)}$.  Suppose, for contradiction,
 $d_j \rho_j = (d_jc_j)/r_{l(j)} \in  (1/r_{l(j-1)})\Z$ where $0
< d_j < s_j$.  Thus, $r_{l(j)} \mid d_jc_jr_{l(j-1)}$. Now, $s_j =
r_{l(j)}/r_{l(j-1)}$, and so $s_j \mid d_jc_j$.  By Lemma
\ref{lem:rhoalpha}, $\gcd(c_j,s_j) = 1$, and so $s_j \mid d_j$.
Since $0 < d_j < s_j$, we have a contradiction.

Now we proceed to show the inductive step.  Suppose $0 \le d_j<
s_j$ for $1 \le j \le i+1$ and $d_{i+1} \not=0$. We write
$$\sum_{j=1}^{i+1} d_j \rho_j = \left(\sum_{j=1}^{i} d_j \rho_j\right) +
d_{i+1} \rho_{i+1}.$$  By the induction hypothesis,
$\sum_{j=1}^{i} d_j \rho_j \in (1/r_{l(i)}) \Z$. Now,
$d_{i+1}\rho_{i+1} \in (1/r_{l(i+1)}) \Z$, and by the previous
paragraph, $d_{i+1}\rho_{i+1} \not\in (1/r_{l(i)}) \Z$.  Thus
$\sum_{j=1}^{i+1} d_j \rho_j \in (1/{r_{l(i+1)}}) \Z -
(1/{r_{l(i)}}) \Z$.
\end{proof}

\section{Representations of Elements of the Value Monoid}
\label{repsmonoid}

 In this section, we
demonstrate that certain elements of ${\Lambda_{}}$ have a unique
representation as a sum of elements of
$\{1, \rho_1, \rho_2, \dots\}$.  Using these representations, we
 prove that ${\Lambda_{}}$ is
generated by $\{1, \rho_1, \rho_2, \rho_3, \dots\}$.  To accomplish this, we
must factor each element of $k[t,y]$ completely as $f(t,y) = q(t) \prod(y-s_i)$
where $s_i$ lies in the algebraic closure of $k(t)$.

An element of $\kq$ is said to be {\bf Puiseux} if it lies in
$k((t^{-1/m}))$ for some positive integer $m$. Puiseux's Theorem
states that the algebraic closure of the field of Laurent series
$k((t^{-1}))$ in $\kq$ precisely consists of all elements of $\kq$
that are Puiseux.   Using Kedlaya's characterization of the
generalized power series that are algebraic over the Laurent power
series field when $k$ has positive characteristic in \cite{ked},
we have the following characteristic-free generalization of
Puiseux's Theorem.

\begin{theorem}
\label{theorem:puiseux}
Let $w\in \kq$ such that no element of its support is divisible by char $k$.
Then $w$
is algebraic over $k((t^{-1}))$ iff $w$ is Puiseux.
\end{theorem}

  The result below follows directly from
techniques found in \cite{abhyankar} and \cite{duval}.

\begin{prop}\label{prop:finitePuiseux}  Let $w= c_1 t^{m_1/n} + \cdots +
c_s t^{m_s/n}$ be a finite Puiseux expansion with ramification
index $n$ where $m_i \in \Zset^*$, $n \in \Zset^+$, and $c_i \in
k^*$.  If $k$ has positive characteristic, then assume that $n$ is
not divisible by char $k$.  Then the minimal polynomial of $w$
over $k(t)$ is $p(y) = \prod_{i=0}^{n-1}(y-w_i) \in k(t)[y],$
where
\[w_i = c_1 ( \zeta^it^{1/n})^{m_1} + \cdots + c_s
(\zeta^it^{1/n})^{m_s}, \] and $\zeta$ is a primitive $n$th
root of unity in $\overline{k}$.
\end{prop}

 The {\bf ramification index}
of a Puiseux series $w\in \kq$ is the smallest positive integer
$r$ such that $w \in k((t^{-1/r}))$. Given $z_1, z_2 \in \kq$, we
say that $z_1$ and $z_2$ {\bf agree to (finite) order } $m\in\N$
if the first $m$ terms of $z_1$ and $z_2$ are identical, but the
$(m+1)$st terms (if they exist) of $z_1$ and $z_2$ are
different.  If we use Theorem \ref{theorem:puiseux} in place of
Puiseux's Theorem, then Proposition 4.6 of \cite{ms3} can be
strengthened to the following characteristic-free form, where we continue
the assumption that no component of the exponent sequence is divisible by
char $k$ as stated in Convention \ref{conventions}.

\begin{prop}\label{prop:lexpformula}
Let $w$ be a Puiseux
series in $\kq$.    Define $p(y) \in \skl[y]$ to be the minimal polynomial of
$w$ over $\skl$ where $w$ agrees with $z$ to order $m$, and none
of the conjugates of $w$ agree with $z$ to a greater order. If $R$
is the ramification index of $w$,  then
\begin{equation}\label{eq:p(z)}
\lexp(p(z)) = \Big(\frac{{{R}}}{r_m}\Big)\Big[  u_m + \lexp(z-w)
\Big]  \ge
\Big(\frac{{{R}}}{r_m}\Big)\Big[  u_m + e_{m+1} \Big] .
\end{equation}
\end{prop}

The simplest polynomials to which we can apply this result are those whose roots
are finite Puiseux series.  We make these calculations explicit in the following lemma.

\begin{lem}\label{lem:rhopos}
If $g(t,y)\in k(t)[y]$ is the minimal polynomial of
\begin{equation*}
 c_1 t^{e_1} + \cdots + c_{l(i)-1}t^{e_{l(i)-1}}\end{equation*}
over $k(t)$, then $\deg_y(g(t,y)) = r_{l(i)-1}$ and $\lexp(g(t,z))
= \rho_i.$
\end{lem}

\begin{proof}
Let $g(t,y) \in k((t^{-1}))[y]$ be the minimal polynomial of
$\sum_{j=1}^{l(i)-1} c_j t^{e_j}$ over $k((t^{-1}))$.
 Since the exponent sequence $\mathbf e$ consists solely of
positive numbers, $g(t,y) \in k[t,y]$ by Proposition
\ref{prop:finitePuiseux}.  Since $\sum_{j=1}^{l(i)-1} c_j t^{e_j}$
has ramification index $r_{l(i)-1}$, it follows from Proposition
\ref{prop:finitePuiseux} that $\deg_y g(t,y) = r_{l(i)-1}$.
Moreover, by Proposition \ref{prop:lexpformula}, $ \lexp(g(t,z))
= \left(\frac{r_{l(i)-1}}{r_{l(i)-1}} \right)(u_{l(i)-1} +
e_{l(i)} ) = \rho_i.$
\end{proof}

We will see that in order to generate ${\Lambda_{}}$, we need only consider images of polynomials
whose roots are finite Puiseux series.  To demonstrate this, we first show that
over the collection of polynomials of a fixed degree in $y$, the polynomials
 that have the smallest image under $\lexp \circ \phi_z$ are those whose roots are finite Puiseux series.

\begin{prop}\label{char0approx}
Let $k$ be a perfect field.
For each nonzero $p(x,y) \in
k[x,y]$, there exists $h(x,y) \in k[x,y]$ such that
the following hold:
\begin{enumerate}
\item[(i)]
 $\deg_y
p(x,y) = \deg_y h(x,y)$,
\item[(ii)] $\lexp(p(t,z)) \ge \lexp(h(t,z))$,
\item[(iii)] the roots of $h(t,y)$ in
$\overline{k((t^{-1}))}[y]$ are finite Puiseux series of the form
$\sum_{j=1}^{l(i)-1} c_j t^{e_j}$.
\end{enumerate}
\end{prop}

\begin{proof}
 First, factor $p(t,y)$ as a polynomial in
$y$ as $p(t,y) =  q(t) \prod_{i=1}^m p_i(t,y),$ where $q(t) \in
k[t]$ and $p_i(t,y)$ is a monic, irreducible element of $\skl[y]$. We
will find $h_i(x,y) \in k[x,y]$ such that $\deg_y p_i(x,y) =
\deg_y h_i(x,y)$, $\lexp(p_i(t,z)) = \lexp(h_i(t,z))$, and the
roots of $h_i(t,y)$ are finite
Puiseux series of the desired form.  It then follows that $h(x,y)
= q(x) \prod_{i=1}^m h_i(x,y)$ satisfies the conditions of the
proposition.

Since $p_i(t,y)$ is
a monic, irreducible element of $\skl[y]$, it is the minimal polynomial of
some generalized power series $\beta \in \kq$.  If $k$ is a field
of characteristic zero, by Puiseux's Theorem (Theorem
\ref{theorem:puiseux}), $\beta$ is Puiseux.  If $k$ has positive
characteristic, $\beta$ is not necessarily Puiseux and the algebraic closure
of $k((t^{-1}))$ is described by Kedlaya in \cite{ked}.
We prove the result by considering two cases:
\begin{enumerate}
\item[Case 1:] No element of $\supp(z)$ is divisible by char $k$.
\item[Case 2:] Some element of $\supp(z)$ is divisible by char $k$.
\end{enumerate}

{\bf Case 1:}
Without loss of generality, we
assume that no conjugate of $\beta$ agrees with $z$ to a higher
order.  We denote this order by $m$, and denote the
ramification index of $\beta$ by $R$, in
which case $r_{m} \mid R$.  As shown in \cite{enum}, $p_i(t,y) \in \skl[y]$ must be a
polynomial of degree $R$.

Let $L$ be the largest index such that $r_{L} = r_{m}$, in which case
 $r_{L+1} > r_{L}$, and so $L+1$ is of the
form $l(\kappa)$ for some $\kappa \in \N$.  Let
$g(t,y) \in
k[t,y]$ be the minimal polynomial of $\sum_{j=1}^{l(\kappa)-1} c_j t^{e_i}$ over $k(t)$. Then by
Lemma \ref{lem:rhopos}, $\deg_y(g(t,y)) = r_{l(\kappa)-1} =
r_{L} = r_{m}$ and $\lexp(g(t,z)) = \rho_{\kappa}.$
Therefore, if we define $h(x,y) = g(x,y)^{(R/r_{m})}$,
then $\lexp(h(t,z)) =  ({R}/{r_{m}}) \rho_{\kappa}$ and
$ \deg_y (h(x,y)) =
({R}/{r_{m}}) \deg_y(g) =
 R =
\deg_y (p_i(x,y)).
$

 Since $r_{L} = r_{m}$, we
know by Lemma \ref{lem:staticuandr} that $u_{L} = u_{m}$.
Moreover, $L  \ge m $, and so $e_{m+1} \ge e_{L+1}$. Thus
by Proposition \ref{prop:lexpformula}, $\lexp(p_i(t,z)) \ge
({R}/{r_{m}})[u_{m} + e_{m+1}] \ge
({R}/{r_{m}})  [u_{L} + e_{L+1}] =
({R}/{r_{m}})  [u_{l(\kappa)-1} + e_{l(\kappa)}] =
({R}/{r_{m}}) \rho_{\kappa} = \lexp(h(t,z)).$

{\bf Case 2 :}
Let char $k = p$.
Let $E$ be the normal closure of $k((t^{-1}))(\beta)/k((t^{-1}))$.
As in the proof of Corollary 9 of \cite{ked}, if $M$ is the integral
closure of $k$ in $E$, then $E$ can be expressed as a tower of Artin-Schreier
extensions over $M((t^{-1/mq}))$, where $q$ is the degree of inseparability
of $E/k((t^{-1}))$.  Since $E$ is normal over $k((t^{-1}))$, and hence over
$k((t^{-1/mq}))$,  the normal closure of $k((t^{-1}))$ must be contained in $E$.
The field $k(\zeta_m)((t^{-1/mq}))= k(\zeta_m, t^{-1/qm})((t^{-1}))$
is the normal closure of $k((t^{-1/mq}))$ (it is the splitting field
of $X^{mq}-t^{-1}$ over $k((t^{-1}))$), and so we have the following
normal extensions:
$$k((t^{-1})) \subset k(\zeta_m)((t^{-1/mq})) \subset E.$$

Define $F=k(\zeta_m)((t^{-1/mq}))$, and let $\tau_\ell \in \gal(F/k((t^{-1})))$ be given by
$t^{1/qm} \mapsto \zeta_m^{q\ell} t^{1/qm}$.   Note that
as $\zeta_m^0, \zeta_m, \dots, \zeta_m^{m-1}$ runs through
all the $m$th roots of unity, so does the list
 $\zeta_m^0, \zeta_m^q, \dots, \zeta_m^{(m-1)q}$ since $\gcd(m,q) = 1$.
Each element
of $\gal(F/k((t^{-1})))$ can be written as $\tau_\ell \mu$ where $\mu \in
\gal(k(\zeta_m)/k)$.  We write the collection of all elements of $\gal(F/k((t^{-1})))$ as
 $\{ \psi_1, \dots, \psi_b\}$.

  Define a homomorphism $\lambda_\ell: \Q \to \overline{k}^*$ by
$\lambda_\ell(ap^n/b) = \zeta_b^{a\ell s}$ where $a\in\Z$,
$b\in\N^*$, $p \nmid  ab$ and $s \equiv p^n \mod b$ (or, if $n<0$,
we require $sp^{-n} \equiv 1 \mod b$). It is straightforward to
show that if $\lambda: \Q \to \overline{k}^*$ is a homomorphism
whose kernel contains $\Z$ and $\mu \in \gal(k(\zeta_m)/k)$, then
\begin{equation}\label{eq:autocoeff}
\sum_{i\in I} x_i t^i \mapsto \sum_{i\in I} \lambda(i) \mu(x_i) t^i
\end{equation}
is a $\overline{k}((t^{-1}))$-automorphism of $\kcq$ (where $I$ is any Noetherian subset of $\Q$).
Given $\psi_j \in \gal(F/k((t^{-1})))$, we write $\psi_j = \tau_\ell \mu$ for some $1 \le \ell \le m$
and $\mu \in
\gal(k(\zeta_m)/k)$.
In case $\lambda = \lambda_\ell$, note that the function in (\ref{eq:autocoeff}) is
an extension of $\psi_j$ to $\kcq$.  We denote the restriction of this function to
$E$ by $\phi_j$.  We will show that $\phi_j$ sends $\overline{k((t^{-1}))}$ to itself,
and since $E$ is a normal extension of $k((t^{-1}))$,
it follows that $\phi_j \in \gal(E/k((t^{-1})))$ is an extension of $\psi_j$.

To show that $\phi_j$ sends $\overline{k((t^{-1}))}$ to itself, we appeal to Kedlaya's
description of the algebraic closure in Corollary 9 of \cite{ked}.  First, we review a few
key ideas from that paper.   The support of any algebraic series must be a set of the form
\begin{equation*}\label{kedsupport}
S_{m,v,c} = \{(1/m)( w+b_1p^{-1} + \cdots + b_{j-1}p^{-j+1} + p^{-n}(b_j p^{-j} + \cdots ))
\mid w \le v, \sum b_i \le c \}
\end{equation*}
where $m \in \N, v,c \ge 0$.  Note that $S_{a,b,c}$
is defined differently than the form given by Kedlaya
since our support is Noetherian rather than well-ordered.
We say that a sequence $c_n$ satisfies a linearized recurrence relation (LRR)
if for some $d_0, \dots, d_k$, for all $n\in \N$,
\begin{equation*}\label{eq:lrr}
d_0 c_n + d_1 c_{n+1}^p + \cdots + d_k c_{n+k}^{p^k} = 0.
\end{equation*}
Let $\sum x_i t^i$ be a series with support $S_{m,v,c}$.  We say $\sum x_i t^i$ is {\bf twist-recurrent} if
for each $w \le v$, $\sum b_i \le c$, the sequence
 $c_n = x_{(1/m)( w+b_1p^{-1} + \cdots + b_{j-1}p^{-j+1} + p^{-n}(b_j p^{-j} + \cdots ))}$
satisfies an LRR.
According to \cite{ked},
the algebraic closure of $k((t^{-1}))$ consists of all twist-recurrent
series $x = \sum x_i t^i$ such that the $x_i$ lie in a finite extension of $k$.

Now suppose $\sum x_i t^i$ is a twist-recurrent series.  We will show
that $\phi_j \left( \sum x_i t^i \right)$ is also twist-recurrent, and so by
the previous paragraph, $\phi_j$ sends $\overline{k((t^{-1}))}$ to itself.
Since $\sum x_i t^i$ is twist-recurrent, it follows that
 $c_n = x_{(1/m)( w+b_1p^{-1} + \cdots + b_{j-1}p^{-j+1} + p^{-n}(b_j p^{-j} + \cdots ))}$
satisfies an LRR of the form
$d_0 c_n + d_1 c_{n+1}^p + \cdots + d_k c_{n+k}^{p^k} = 0$.
To show that $\phi_j \left( \sum x_i t^i \right)$ is twist-recurrent,
we must prove that $\lambda(f(n)) \mu (c_n)$ satisfies an LRR where
$f(n) = (1/m)( w+b_1p^{-1} + \cdots + b_{j-1}p^{-j+1} + p^{-n}(b_j p^{-j} + \cdots ))$,
$\lambda = \lambda_\ell$ for some $\ell$ and $\mu \in \gal(k(\zeta_m)/k)$.
If $c_n$ satisfies the LRR $\sum_{i=0}^k d_i c_{n+i}^{p^i} =0$,
it follows that $0 = \mu \left(\sum_{i=0}^k d_i c_{n+i}^{p^i} \right)
= \sum_{i=0}^k \mu(d_i) \mu(c_{n+i})^{p^i}$, and so $\mu (c_n)$ satisfies
an LRR.  Thus we only have to show that if $c_n$ satisfies an LRR,
then so does $c_n' = \lambda_\ell(f(n)) c_n$.

Now suppose $c_n$ satisfies the LRR $\sum_{i=0}^k d_i c_{n+i}^{p^i}= 0$.
Rewrite $w+ b_1p^{-1} + \cdots + b_{j-1}p^{-j+1}$ as $\frac{\alpha_1}{p^{m_1}}$ where
$p \nmid \alpha_1$ and $m_1 \le j-1$.
If we rewrite $b_j p^{-j} + b_{j+1} p^{-j-1}\cdots $ as
$\frac{\alpha_2}{p^{m_2}}$ where
$p \nmid \alpha_2$ and $m_2 \ge j$, then
\begin{equation*}
f(n) = \frac{\alpha_1 p^{m_2+n} + \alpha_2 p^{m_1}}{mp^{m_1+m_2+n}}
= \frac{\alpha_1 p^{m_2-m_1+n} + \alpha_2}{mp^{m_2+n}}.
\end{equation*}
If we define $s_n, d_1, d_2$ so that $s_np^n \equiv 1 \mod m$,
$d_1p^{m_1} \equiv 1 \mod m$, and $d_2p^{m_2} \equiv 1 \mod m$, then
$\lambda_\ell(f(n)) =\zeta_m^{(\alpha_1p^{m_2-m_1+n}+\alpha_2)s_nd_2}
= \zeta_m^{\alpha_1d_1} \cdot \zeta_m^{\alpha_2d_2s_n}$,
and so if we define $d_i' = \zeta_m^{-\alpha_1d_1p^i}d_i$, then
$$\sum_{i=0}^k d_i' c_{n+i}'^{p^i} =
\sum_{i=0}^k (\zeta_m^{-\alpha_1d_1p^i})d_i({\zeta_m^{\alpha_1d_1} \cdot \zeta_m^{\alpha_2d_2s_{n+i}}})^{p^i} c_{n+i}^{p^i} =
\sum_{i=0}^k (\zeta_m^{-\alpha_1d_1p^i})d_i({\zeta_m^{\alpha_1d_1p^i})(\zeta_m^{\alpha_2d_2s_{n+i}}})^{p^i} c_{n+i}^{p^i},$$
which simplifies as
$$
\sum_{i=0}^k (\zeta_m^{\alpha_2d_2s_{n+i}})^{p^i}d_i c_{n+i}^{p^i} =
\sum_{i=0}^k (\zeta_m^{\alpha_2d_2s_{n}{s_i}})^{p^i}d_i c_{n+i}^{p^i} =
(\zeta_m^{\alpha_2d_2s_{n}}) \sum_{i=0}^k d_i c_{n+i}^{p^i} = 0,
$$
and so $c_n'$ satisfies an LRR.

So far,  we have shown that $\phi_j$ sends $\overline{k((t^{-1}))}$ to itself,
and since $E$ is a normal extension of $k((t^{-1}))$,
we know $\phi_j \in \gal(E/k((t^{-1})))$ is an extension of $\psi_j$.
  Let $\{\sigma_1, \dots, \sigma_d\}$ be the
complete collection of $F$-automorphisms of
$E$.
Since $E/F$ and $F/k((t^{-1}))$ are normal extensions, a routine
exercise shows that
 the collection
$\{\phi_i\sigma_j \mid 1 \le i \le b, 1 \le j \le d \}$ consists of all $k((t^{-1}))$-automorphisms of $E$.
 Since $q$ is the degree of inseparability of
$E$ over $k((t^{-1}))$,  the minimal polynomial
$m_\beta$ of $\beta$ over $k((t^{-1}))$ can be factored as
\begin{equation*}\label{eq:minpoly}
m_\beta(t,y) = \prod_{i=1}^{d}
\left( \prod_{j=1}^{b} (y -  \phi_j \sigma_i\beta)
\right)^{q}.
\end{equation*}  For any series $s = \sum_{i \in I} c_i t^{e_i}$, we define an
associated Puiseux series by $\mathcal{P}(s) = \sum_{i \in J}
c_i t^{e_i}$ where $J = \{ a/b \in I \mid a\in \Z,b \in \N^*
\mbox{ and } p\nmid b \}$ and remainder by $\mathcal{R}(s) = s
- \mathcal{P}(s)$.
 Since no
component of the ramification sequence of $z$ is divisible by $p$, we
obtain
\begin{equation}\label{eq:boundle}
\lexp(z-\phi_j \sigma_i \beta)
= \lexp(z- {\mathcal P}(\phi_j \sigma_i \beta) -
{\mathcal R}(\phi_j \sigma_i \beta))
\ge \lexp(z-{\mathcal P}(\phi_j \sigma_i \beta)).
\end{equation}
Since $\phi_j$ is of the form (\ref{eq:autocoeff}),
for any series $s \in\kq$,
$\mathcal{P}(\phi_js) = \phi_j({\mathcal P}(s))$.  Applying this to
(\ref{eq:boundle}), we obtain
\begin{equation*}\label{eq:boundle2}
\lexp(z-\phi_j \sigma_i \beta)
\ge \lexp(z-\phi_j {\mathcal P}(\sigma_i \beta)).
\end{equation*}
Of all the conjugates $\phi_j({\mathcal P}(\sigma_i \beta))$ of
${\mathcal P}(\sigma_i \beta)$ over $F$, choose
$\alpha_i$ to be the one that agrees with $z$ to the highest order.
Note that $\prod_{j=1}^b (y-\phi_j \alpha_i)$ must be
of the form $m_{\alpha_i}(t,y)^{\ell_i}$  where $m_{\alpha_i}(t,y)$ is
the minimal polynomial of $\alpha_i$ over $k((t^{-1}))$ and $\ell_i \in \N$.
Since $\alpha_i$ is a Puiseux series such that no element of its support
is divisible by $p$, we have reduced the problem to Case 1, and the proof is
complete.
\end{proof}

Now, we define a
sequence of rational numbers that give the minimal possible value
of an image of a polynomial of degree $d$ under the map $\lexp
\circ \varphi_z$.

\begin{definition}\label{def:lambda}
For each natural number $d$,
\begin{equation*}
\lambda_d   :=    \min \{ \lexp(f(t,z)) \mid f\in k[x,y]^* \mbox{
\rm and } \deg_y(f(x,y)) = d \}. \end{equation*}
\end{definition}

\begin{lem}\label{lem:lambdad}
Let $k$ be a perfect field.
For any positive integer $d$,
 \begin{equation}\lambda_d =
\lexp\left(\prod_{j=1}^w f_j(t,z)^{d_j} \right)
\end{equation} where
$w$ is a positive integer,  the exponent $d_j$ is nonnegative, and
$f_j$ is the minimal polynomial of $\sum_{i=1}^{{l(j)-1}} c_i
t^{e_i}$ over $k(t)$.  Moreover, $d= \sum d_j \deg_y(f_j(x,y))$.
\end{lem}

\begin{proof}
 By
the definition of $\lambda_d$, there exists $p(x,y) \in k[x,y]$
such that $\deg_y(p(x,y))= d$ and $\lexp(p(t,z)) = \lambda_d.$ By
Proposition \ref{char0approx}, there
exists $h(x,y)$ such that $\lambda_d = \lexp(p(t,z)) \ge
\lexp(h(t,z))$, $\deg_y(h(x,y)) = d$,
and $h(t,y)$ has finite Puiseux series as roots.  Thus, by the definition
of $\lambda_d$, $ \lambda_d = \lexp(h(t,z)).$
 Since $h(x,y)$ is a
product of minimal polynomials of finite Puiseux series, we can
write $h$ as
$h(t,z) = \prod_{j=1}^w f_j(t,z)^{d_j},
$ where $w$ is a positive
integer, and for each $1 \le j \le w$, the exponent $d_j$ is
nonnegative, and $f_j$ is the minimal polynomial of
$\sum_{i=1}^{{l(j)-1}} c_i t^{e_i}$ over $k(t)$.
\end{proof}

Using this lemma, we can produce a unique representation for each $\lambda_d$ in terms of the
monoid generating sequence.

\begin{prop}\label{lem:minrep}
Let $k$ be a perfect field.
For any positive integer $d$, $\lambda_d$ can be uniquely expressed in the form
\begin{equation}\label{eq:lambdaunique}
\lambda_d = \sum_{j=1}^w d_j \rho_j,
\end{equation}
where $w$ is a positive integer, and for each $1 \le j \le w$, we
have
\begin{equation}\label{eq:boundd}
0 \le d_j < s_j.
\end{equation}
In this case,
\begin{equation*}\label{eq:decomposed}
d = \sum_{j=1}^d d_j r_{l(j-1)}.
\end{equation*}
\end{prop}

\begin{proof}
By
Lemma \ref{lem:lambdad}, there
exists $h(x,y) \in k[x,y]$ such that $\lambda_d =
\lexp(h(t,z))$, $\deg_y(h(x,y)) = d$, and
\begin{equation*}\label{eq:hprod}
h(t,z) = \prod_{j=1}^w f_j(t,z)^{d_j},
\end{equation*}
 where $w$ is a positive
integer, and for each $1 \le j \le w$, the exponent $d_j$ is
nonnegative, and $f_j$ is the minimal polynomial of
$\sum_{i=1}^{{l(j)-1}} c_i t^{e_i}$ over $k(t)$.
By Lemma \ref{lem:rhopos}, $\deg_y f_j(x,y) = r_{l(j)-1}$ and $\lexp(
f_j(t,z)) = \rho_j$, and so
\begin{equation*}\label{eq:lambdanew2}
\lambda_d  = \lexp\left(\prod_{j=1}^w f_j(t,z)^{d_j} \right) =
\sum_{i=1}^w d_j \lexp(f_j(t,z)) = \sum_{i=1}^w d_j \rho_j
\end{equation*}
and
\begin{equation*}
d = \deg_y h(x,y)  = \sum_{j=1}^w d_j \deg_y f_j(x,y) =
\sum_{j=1}^w d_j r_{l(j)-1} = \sum_{j=1}^w d_j r_{l(j-1)}.
\end{equation*}

Next we show that each $d_j$ satisfies the bounds given by
(\ref{eq:boundd}). Suppose for contradiction, for some $k$, $d_k
\ge s_k = r_{l(k)}/r_{l(k-1)}$. Define
$$D_j \ = \ \left\{ \begin{array}
{l@{\quad}cl}
d_j + 1  & & \mbox{if } j= k+1; \\
d_j - s_j & & \mbox{if } j= k; \\
d_j & & \mbox{otherwise}.
\end{array} \right. $$
Using this in conjunction with the recurrence relation given in
Lemma \ref{lem:therecurrence}, we obtain
\begin{eqnarray*}
\sum_{j=1}^w d_j \rho_j - \sum_{j=1}^w D_j \rho_j & = & (d_k
- D_k) \rho_k  +(d_{k+1} - D_{k+1})  \rho_{k+1} \\
& = &  s_k \rho_k - \rho_{k+1} \\
& = & e_{l(k)} - e_{l(k+1)},
\end{eqnarray*}
and so
\begin{eqnarray*}
\sum_{j=1}^w d_j r_{l(j-1)}  - \sum_{j=1}^w D_j r_{l(j-1)} & = &
(d_k
- D_k) r_{l(k-1)} + (d_{k+1} - D_{k+1})  r_{l(k)} \\
& = & s_k r_{l(k-1)} - r_{l(k)} \\
& = & 0.
\end{eqnarray*}
These equations in conjunction with Lemma \ref{lem:rhopos} yield
$$\lexp\left(\prod_{j=1}^w f_j(t,z)^{D_j} \right)
= \sum_{j=1}^w D_j \rho_j = \sum_{j=1}^w d_j \rho_j - e_{l(k)} +
e_{l(k+1)} < \sum_{j=1}^w d_j \rho_j = \lexp(h)
$$
and
$$\deg\left(\prod_{j=1}^w f_j(t,z)^{D_j}\right)
= \sum_{j=1}^w D_j \deg(f_j) = \sum_{j=1}^w D_j r_{l(j-1)} =
\sum_{j=1}^w d_j r_{l(j-1)}  = \deg(h).$$ However, $\lexp(h) =
\lambda_d$, and so we have contradicted the minimality of
$\lexp(h)$. Thus $0 \le d_j < s_j$ for each $1 \le j \le w$, and
so we have proved the bounds given by (\ref{eq:boundd}).

Finally, we demonstrate that the expression for $\lambda_d$ in
(\ref{eq:lambdaunique}) is uniquely determined. Suppose we are
given two representations for $\lambda_d$:
\begin{equation*}
\lambda_d = \sum_{j=1}^w d_j \rho_j = \sum_{j=1}^w d_j' \rho_j
\end{equation*}
where $0 \le d_j, d_j' < s_j$. If we define $\Delta_j = d_j -
d_j'$, then $\sum_{j=1}^w \Delta_j \rho_j = 0$ and $|\Delta_j| <
s_j$. Multiply the expression by $r_{l(w-1)}$, and we see
$$\left(\sum_{j=1}^{w-1} r_{l(w-1)} \Delta_j\rho_j\right)
+ r_{l(w-1)} \Delta_w \rho_w =0.$$ However, $r_{l(w-1)} \Delta_j
\rho_j \in \Z$ for $j \le w-1$, and so $r_{l(w-1)} \Delta_w \rho_w
\in \Z$.  Now write $\rho_w$ as $c_w/r_{l(w)}$ where $c_w \in \N$.
Then $r_{l(w-1)} \Delta_w c_w/r_{l(w)} \in \Z$, and so $s_w =
\frac{r_{l(w)}}{r_{l(w-1)}} \mid \Delta_wc_w$. Since $s_w$ and
$c_w$ are relatively prime by Lemma \ref{lem:rhoalpha}, $s_w \mid
\Delta_w$. However, $|\Delta_w| < s_w$, and so $\Delta_w = 0$.
Thus, $\sum_{j=1}^{w-1} \Delta_j \rho_j = 0$. Repeating this
argument, we find $\Delta_{w-1} = \Delta_{w-2} = \cdots = \Delta_1
= 0$, and so $d_j = d_j'$ for all $1 \le j \le w$.
\end{proof}

The idea that each $\lambda_d$ has a unique representation can be
extended further.  In fact, there is a natural bijective
correspondence between representations of natural numbers and
representations of terms of the form $\lambda_d$.
First, we state the following simple lemma without
proof.

\begin{lem}\label{lem:baserep}
Let $b_0, b_1, b_2, b_3, \dots$ be a sequence of positive integers
such that $b_0 = 1, b_{i+1} > b_i$ and $b_i \mid b_{i+1}$ for all
$i$. Then every positive integer $n \in \N$ has a unique
representation of the form
$$d = \sum_{i=0}^w d_i b_i,$$
where $w$ is a positive integer, $d_w \not=0$, and $0 \le d_i <
b_{i+1}/b_{i}$.
\end{lem}

For example, if $b_i = {10}^i$, then this says that every positive
integer has a unique base 10 representation.  Using this lemma, we produce a method
for quickly computing $\lambda_d$.

\begin{prop}\label{prop:correspondence}
\label{eq:lambdad}
Let $k$ be a perfect field.
Given a positive integer $w$ and $0 \le d_j <
s_j$ for each $1 \le j \le w$,
$$d = \sum_{j=1}^w d_{j}
r_{l(j-1)}
 \ \Leftrightarrow \ \lambda_d = \sum_{j=1}^w d_j \rho_j.$$
\end{prop}

\begin{proof}
The reverse implication follows directly from Proposition
\ref{lem:minrep}.  For the forward implication, suppose we are
given $d= \sum_{j=1}^w d_j r_{l(j-1)}$ where $0 \le d_j < s_j$. By
Proposition \ref{lem:minrep}, $\lambda_d$ is of the form
$\lambda_d = \sum_{j=1}^{w'}d_j' \rho_j$  where
$d=\sum_{j=1}^{w'}d_j' r_{l(j-1)}$. By the uniqueness promised by
Lemma \ref{lem:baserep}, $w=w'$ and $d_j=d_j'$ for all $1\le j \le
w$. Thus $\lambda_d = \sum_{j=1}^w d_j \rho_j$.
\end{proof}

\section{Construction of the Value Monoid} \label{construction}

The goal of this section is to  describe the value
monoid ${\Lambda_{}}$ explicitly in terms of the sequences $\{\lambda_i\}_{i\in\N}$
and $\{\rho_i\}_{i\in\N}$.  Throughout the remainder, in addition to Convention \ref{conventions}, we assume that $k$ is
a perfect field and $\{\lambda_i\}_{i\in\N}$ is given by Definition \ref{def:lambda}.  We begin by showing that
$\{\lambda_i\}_{i\in\N}$ is an increasing sequence.

\begin{lem}\label{lem:lambdainc}
The sequence $\lambda_0, \lambda_1, \lambda_2, \dots$ is
increasing.
\end{lem}

\begin{proof}
 We will show that $\lambda_{d+1} >
\lambda_d$ for all $d$. By Proposition \ref{lem:minrep}, we can
write $\lambda_d = \sum_{j=1}^w d_j \rho_j$ where
 $0 \le d_j < s_j$ and
\begin{equation*}
d = \sum_{j=1}^w d_j r_{l(j-1)}.
\end{equation*}
We now consider different cases, depending on the size of the
coefficients $d_j$.

\vskip 0.3 cm \noindent {\bf Case 1: \ } First we consider the
case $d_j = s_j-1$ for all $j$.  Then $d=\sum_{j=1}^w (s_j-1)
r_{l(j-1)}$, and so by Lemma \ref{lem:sumram}, $d+1 = r_{l(w)}$.
Thus by Proposition \ref{prop:correspondence}, $\lambda_{d+1} =
\rho_{w+1}$ and $\lambda_d = \sum_{j=1}^w d_j \rho_j$, and so by Lemma \ref{lem:difflams}, $\lambda_{d+1} -
\lambda_d  =  \rho_{w+1} -\sum_{j=1}^w (s_j-1) \rho_{j} =
e_{l(w+1)} > 0.$

\vskip 0.3 cm \noindent {\bf Case 2: \ }  Consider the case $d_1 <
s_1 -1$.  Now $d+1 = (d_1+1)r_{l(0)}+ \sum_{j=2}^w d_j
r_{l(j-1)}$, and so by Proposition \ref{prop:correspondence},
$\lambda_{d+1} = (d_1+1)\rho_1 + \sum_{j=2}^w d_j \rho_j$. Thus
$\lambda_{d+1}-\lambda_d = (d_1+1)\rho_1 - d_1\rho_1 = \rho_1 >
0$.

\vskip 0.3 cm \noindent {\bf Case 3: \ } Finally we consider the
case where there exists an index $v>1$ such that $d_v < s_v-1$ and
for $j<v$, $d_j = s_j-1$.   Write $\lambda_d$ as $\lambda_d =
\sum_{j=1}^{v-1} (s_j-1) \rho_j +  \sum_{j=v}^{w} d_j\rho_j.$
 By Proposition \ref{prop:correspondence},
 $d =
\sum_{j=1}^{v-1} (s_j-1) r_{l(j-1)} + \sum_{j=v}^w d_j
r_{l(j-1)},$
 and so by Lemma \ref{lem:sumram},
$$d+1 = 1 + \sum_{j=1}^{v-1} (s_j-1) r_{l(j-1)} + \sum_{j=v}^w d_j
r_{l(j-1)}= r_{l(v-1)} + \sum_{j=v}^w d_j r_{l(j-1)}
= (d_v +1) r_{l(v-1)} + \sum_{j=v+1}^w d_j r_{l(j-1)}.$$
 Therefore, by Proposition \ref{prop:correspondence},
$ \lambda_{d+1} =  (d_v+1) \rho_{v} + \sum_{j=v+1}^w d_j \rho_j,$ and so
$\lambda_{d+1} - \lambda_d = (d_v + 1) \rho_v +  \sum_{j=v+1}^w d_j \rho_j - (\sum_{j=1}^{v-1}(s_j-1) \rho_j +
 \sum_{j=v}^w d_j \rho_j)=  \rho_v
- \sum_{j=1}^{v-1} (s_j-1) \rho_j $.
By
Lemma
\ref{lem:difflams},  this is simply $e_l(v)$, which is positive.
\end{proof}

Given a submonoid $M$ of a commutative monoid $N$, we define an
equivalence relation on $N$ by setting $n_1 \sim_M n_2$ if and
only if there exist $m_1, m_2 \in M$ such that $m_1+n_1 =
m_2+n_2$. Denote by $N/M$ the collection of all equivalence
classes under this relation, and  define a quotient map $\pi$ from
$N$ to $N/M$ that sends $n$ to the equivalence class containing
$n$. The set $N/M$ has an additive monoid structure where we
define $\pi(n_1) + \pi(n_2) = \pi(n_1+n_2)$.

Given a polynomial $f(x,y) \in k[x,y]$, we define
${{\mbox{deg}}}_y(f(x,y))$ to be the smallest $d\ge0$ such that
$f(x,y)\in k[x]y^d+k[x]y^{d-1}+\cdots+k[x]y+k[x]$, and we denote
\begin{equation}\label{def:BigLambda}
\Lambda_d(z) = \{ \lexp(f(t,z)) \mid f\in k[x,y]^* \mbox{ \rm and
} \deg_y(f(x,y)) \le d \}.\end{equation}

Using this notation, we show that any pair of terms of the
sequence $\{\lambda_i\}_{i\in\N}$
are inequivalent modulo $\Z$.

\begin{prop}\label{prop:minrepsinequiv}
For all $i \not= k$, $\lambda_i \not\sim_\Z \lambda_k$.
\end{prop}

\begin{proof}
Suppose $\lambda_i \sim_\Z \lambda_k$.
 By Proposition \ref{lem:minrep}, for some positive integer $w$ we can write $\lambda_i =
\sum_{j=1}^w d_j \rho_j$ and $\lambda_k = \sum_{j=1}^w d_j'
\rho_j$ where $0 \le d_j, d_j' < s_j$.  For each $1 \le j \le w$,
we write $\rho_j = c_j/r_{l(j)}$, where $c_j$ and $s_j$ are
relatively prime, as promised by Lemma \ref{lem:rhoalpha}.

If we define $\Delta_j = d_j - d_j'$, then $|\Delta_j| < s_j =
r_{l(j)}/r_{l(j-1)}$ and $\lambda_i - \lambda_k = \sum_{j=1}^w
\Delta_j \rho_j \sim_\Z 0$. Multiply the expression by
$r_{l(w-1)}$ to obtain
\begin{eqnarray}\label{eq:congcancel}
\left(\sum_{j=1}^{w-1} r_{l(w-1)} \Delta_j\rho_j\right) +
r_{l(w-1)} \Delta_w \rho_w \sim_\Z 0. \end{eqnarray} However,
$r_{l(w-1)} \Delta_j \rho_j \in \Z$ for $j \le w-1$ since $\rho_j
\in (1/r_{l(j)})\Z$, and so by (\ref{eq:congcancel}), $r_{l(w-1)}
\Delta_w c_w/r_{l(w)} =  r_{l(w-1)} \Delta_w \rho_w \in \Z$. That
is, $\Delta_w c_w /s_w = r_{l(w-1)} \Delta_w c_w/r_{l(w)} \in \Z$,
and so $s_w \mid \Delta_w c_w$. Since $s_w$ and $c_w$ are
relatively prime, $s_w \mid \Delta_w$. However, $|\Delta_w| <
s_w$, and so $\Delta_w = 0$. Thus, $\sum_{j=1}^{w-1} \Delta_i
\rho_j \sim_\Z 0$. Repeating this argument, we find $\Delta_{w-1}
= \Delta_{w-2} = \cdots = \Delta_1 = 0$, and so $\lambda_i =
\lambda_k$. By Lemma \ref{lem:lambdainc}, $i=k$.

\end{proof}

We quote
the following result from \cite{ms2}.

\begin{theorem}\label{theorem:digistight}
For every positive integer $n$, the quotient $\Lambda_d/\Lambda_0$ has cardinality one greater
than that of $\Lambda_{d-1}/\Lambda_0$, or equivalently,
$\Lambda_d/\Lambda_0$ has cardinality $d+1$.
\end{theorem}

Using this theorem in conjunction with Proposition  \ref{prop:minrepsinequiv}, we compute the quotient
$\Lambda_d/\Lambda_0$.

\begin{cor}\label{cor:quotientgen}
The quotient $\Lambda_d / \Lambda_0$ consists precisely of the
images of $\lambda_0, \dots, \lambda_d$.
\end{cor}

\begin{proof}
Since $\lambda_0, \dots, \lambda_d \in \Lambda_d$, we know by
Proposition \ref{prop:minrepsinequiv} that the images of
$\lambda_0, \dots, \lambda_d$ are distinct in $\Lambda_d /
\Lambda_0$. By Theorem \ref{theorem:digistight}, these images
constitute the entire quotient $\Lambda_d / \Lambda_0$.
\end{proof}

For each $m \in {\Lambda_{}}$, we make the following definition:
\begin{equation}\label{def:lambdam}
\lambda(m) = \min \{ r \in {\Lambda_{}} \mid r \sim_\Z m \}.
\end{equation}
The next two results allow us to relate terms of the sequence
$\{\lambda_i\}_{i\in\N}$ with elements in
the image of the map $\lambda: {\Lambda_{}} \to {\Lambda_{}}$.

\begin{prop}\label{prop:minrepcross}
For all $i \in \N$, there exists $m \in {\Lambda_{}}$ such that $\lambda_i
= \lambda(m)$.
\end{prop}

\begin{proof}
We prove the following equivalent statement:
 for all $i \in
\N, m \in {\Lambda_{}}$, if $m \sim_\Z \lambda_i$, then $\lambda_i \le m$.
Let $i \in \N$, $m\in {\Lambda_{}}$ such that $m \sim_\Z \lambda_i$. Let
$j$ be the smallest index such that $m \in \Lambda_j$.  Suppose, for contradiction,
 $j<i$.  Since the image of
$m$ must lie in the quotient $\Lambda_j /\Lambda_0$, by Corollary
\ref{cor:quotientgen} it follows that $m \sim_\Z \lambda_t$ for some $t\le
j <  i$. Thus, $\lambda_i \sim_\Z \lambda_t$, which contradicts
Proposition \ref{prop:minrepsinequiv}.
Therefore,  $j \ge
i$, and so by Lemma \ref{lem:lambdainc}, $m \ge \lambda_j \ge
\lambda_i$.
\end{proof}

\begin{prop}\label{prop:minstrata}
For all $m\in {\Lambda_{}}$, there exists $i \in \N$ such that $\lambda_i =
\lambda(m)$.
\end{prop}

\begin{proof}
Let $m \in {\Lambda_{}}$.   Now $m \in \Lambda_j$ for some $j \in \N$, and so by
Corollary \ref{cor:quotientgen}, $m \sim_\Z \lambda_i$ for some $i
\in \N$. By Proposition \ref{prop:minrepcross}, $\lambda_i =
\lambda(m')$ for some $m' \in {\Lambda_{}}$.  Thus $\lambda_i \sim_\Z m
\sim_\Z m'$, and so $\lambda_i = \lambda(m') = \lambda(m)$.
\end{proof}

We are now in a position to decompose the value monoid as a disjoint
union of cosets of $\N$.

\begin{theorem}\label{thm:monoidreplambda}
If the exponent sequence of $z \in \kq$ is strictly positive, then
the value monoid is the disjoint union
\begin{equation*}\label{mzuplus}
{\Lambda_{}} = \bigcup_{d=0}^\infty (\N + \lambda_d).
\end{equation*}
\end{theorem}

\begin{proof}
Given $m \in {\Lambda_{}}$, there exists an index $d$ such that $\lambda_d
= \lambda(m)$ by Proposition \ref{prop:minstrata}.  Therefore,
$m - \lambda_d \in \N$, and so $m \in \N + \lambda_d$. The reverse
containment follows directly from the fact that $\lambda_d \in
{\Lambda_{}}$.
 The sets are disjoint due to Proposition
\ref{prop:minrepsinequiv}.
\end{proof}

Combining Theorem \ref{thm:monoidreplambda} and Proposition
\ref{lem:minrep}, we obtain the following.

\begin{theorem}\label{Mzunique}
Each element $m \in {\Lambda_{}}$ has a unique representation of the form
\begin{equation}\label{eq:Mzunique}
m = n + \sum_{j=1}^w d_j \rho_j,
\end{equation}
where $n\in \N$ and for each $1 \le j \le w$, $0 \le d_j < s_j.$
\end{theorem}

A weaker form of this theorem was stated earlier as
Theorem \ref{thm:valmonoid}.

\section{Algorithms} \label{algorithms}

In this section, we develop algorithms to make computations involving the value monoid ${\Lambda_{}}$.
It was shown in  \cite{ms3} that ${\Lambda_{}}$ is well-ordered, and so
 $\lexp \circ \varphi_z$ is suitable relative to $\kx$ as described
in Definition \ref{def:suitableval}, and we can use $\lexp \circ \varphi_z$
in the algorithms described in Section \ref{introduction}.
Throughout this section we refer to the composite maps $\lexp \circ \varphi_z$
and $\lc \circ \varphi_z$  as $\lexpz$ and $\lcz$, respectively.

To begin, given a rational number $m \in \Q$, we would like to decide
whether $m\in {\Lambda_{}}$, and in case it is, express it in terms of the
generators $1, \rho_1, \rho_2, \dots$.  To accomplish this, we first prove a lemma.

\begin{definition}
For each $i \in \N$, define
$$\Omega_i = \{ n+ \sum_{j=1}^i d_j \rho_j \mid, n\in \N, 0 \le d_j < s_j\}.$$
\end{definition}

\begin{lem}\label{lem:NorZ}
$${\Lambda_{}} \cap \Z \cdot \{1, \rho_1, \rho_2, \rho_3, \dots, \rho_i\} =
 \Omega_i.$$
\end{lem}

\begin{proof}
The containment `$\supset$' being obvious, we only consider the
case `$\subset$'.  Let $m \in {\Lambda_{}} \cap \Z \cdot \{1, \rho_1,
\rho_2, \rho_3, \dots, \rho_i \}$.  By Theorem
\ref{thm:monoidreplambda}, there is a unique pair $n, d \in \N$
such that $m = n +  \lambda_d$. Thus $\lambda_d \in \Z \cdot \{1,
\rho_1, \rho_2, \dots, \rho_i \}$, and so by Lemma
\ref{cor:rhoresidue},  $\lambda_d \in (1/r_{l(i)}) \Z$.

By Theorem \ref{Mzunique}, there exists a
smallest $k\in \N$
such that  $\lambda_d \in \Omega_k$.
Suppose, for contradiction, that $k > i$. Then by Lemma
\ref{lem:lambdaresidue}, $\lambda_d \in (1/ r_{l(k)}) \Z -  (1/
r_{l(k-1)}) \Z  \subset (1/ r_{l(k)}) \Z - (1/ r_{l(i)})\Z$, which
contradicts our assertion that $\lambda_d \in (1/r_{l(i)}) \Z$.
Therefore, $i=k$, and so $\lambda_d \in \N \cdot \{1, \rho_1,
\dots, \rho_i \}$.
\end{proof}

We have the following corollary.

\begin{cor}\label{omegaclosed}
The set $\Omega_i$ is closed under addition.
\end{cor}

Given a positive rational number $m$, write $m$ as
${a}/{b}$ where $a,b$ are relatively prime positive integers.
If $m \in \N$, then it is automatically in ${\Lambda_{}}$, and so we can
assume that $b > 1$. Our goal is to decide using modular
arithmetic whether it is possible that $m \in {\Lambda_{}}$.  First, find
the smallest $i$ such that $b \mid r_{l(i)}$.  The set of all
$\Z$-linear combinations of $1, \rho_1, \dots, \rho_{i-1}$ is
precisely the set $\frac{1}{r_{l(i-1)}}\Z$.  Since $b$ does not
divide $r_{l(i-1)}$, it cannot possibly be an $\N$-linear
combination of $1, \rho_1, \dots, \rho_{i-1}$.  Now suppose $m$ is
a $\Z$-linear combination of $1, \rho_1, \dots, \rho_j$ where $j>
i$.  However, since $b \mid r_{l(i)}$, it follows that $m \in
(1/r_{l(i)}) \Z = \Z \cdot \{1, \rho_1, \dots, \rho_i\}$.
If $m \in \Lambda$, then by Lemma \ref{lem:NorZ}, there exist $n, d_1,
\dots, d_i \in \N$ such that
\begin{equation*}
m = n + \sum_{j=1}^i d_j \rho_j
\end{equation*}
where $0 \le d_j < s_j$ for $1 \le j \le i$ and $d_i\not=0$.
From this discussion, we have the following algorithm.

\begin{alg}\label{alg:decomposem}
Let $m$ be a positive rational number.  The following algorithm
determines whether $m \in {\Lambda_{}}$.  If $m \in {\Lambda_{}}$, then the
algorithm produces a decomposition of $m$ as a linear combination
of $1, \rho_1, \dots, \rho_i$. Set $\rho_i = c_i/r_{l(i)}$.

\begin{enumerate}
\item[(1)] Write $m$ as $a/b$ where $a,b$ are relatively prime,
positive integers.

\item[(2)] Define $i$ to be the smallest index such that $b \mid
r_{l(i)}$.

\item[(3)] Define $m^{(i)} = m$.

\item[(4)]  Try to solve the congruence $c_id_i \equiv r_{l(i)}
\mod s_i$ for $d_i$ where $0 \le d_i < s_i$.  If there are no
solutions, then $m\not\in {\Lambda_{}}$.

\item[(5)]  For $j = i-1, i-2, \dots, 1$, define $m^{(j)} = m^{(j+1)} -
d_{j+1}\rho_{j+1}$ and try to solve the congruence $c_jd_j \equiv
r_{l(j)} \mod s_j$ for $d_j$ where $0 \le d_j < s_j$.  If any of
the congruences fail to yield a solution, then $m\not\in {\Lambda_{}}$.

\item[(6)] Define $n = m^{(1)} - d_1 \rho_1$. Then
$m = n + \sum_{j=1}^i d_j \rho_j.$
If $n \not\in \N$, then $m\not\in {\Lambda_{}}$.  If $n\in {\Lambda_{}}$, then we
have a decomposition of the desired form.
\end{enumerate}
\end{alg}

Once we have a test for whether a rational number is in the value monoid, we need
to be able to determine one of its preimages under the valuation.  The following algorithm
accomplishes this task.

\begin{alg}\label{alg:monoidtopoly}
Let $m \in {\Lambda_{}}$.  This algorithm constructs $p(x,y) \in k[x,y]$
such that $\lexpz(p(x,y)) = m$.

\begin{enumerate}
\item[(1)] Using Algorithm \ref{alg:decomposem}, write $m =
n + \sum_{j=1}^i d_j \rho_j.$
\item[(2)] For each $1 \le j \le i$, use Proposition \ref{prop:finitePuiseux}
to compute $p_j(x,y)$,  the minimal polynomial of
$\sum_{j=1}^{l(i)-1} c_j x^{e_j}$ over $k(x,y)$.
\item[(3)] Define $p(x,y) = x^n \prod_{j=1}^i p_j(x,y)^{d_j}$.  By Lemma \ref{lem:rhopos},
$\lexpz(p(x,y)) = m$.
\end{enumerate}
\end{alg}

The following algorithm describes how to
perform division in $k[x,y]$ relative to $\lexpz$.

\begin{alg}\label{alg:division}
Let $f,g \in \kx$.  This algorithm constructs $h \in k[x,y]$
such that $\lexpz(f-gh) < \lexpz(f)$ provided that such an $h$
exists.

\begin{enumerate}
\item[(1)]  Compute $m = \lexpz(f) - \lexpz(g)$.
\item[(2)]  Use Algorithm \ref{alg:decomposem} to determine whether $m \in \Lambda$.
If $m\not\in \Lambda$, then  $h$ does not exist.
\item[(2)] Using Algorithm \ref{alg:monoidtopoly}, find $p(x,y) \in k[x,y]$ such that $\lexpz(p) = m$.
\item[(3)] Define $h(x,y) = (\lcz(f)/\lcz(gp)) p(x,y)$.  Then $\lcz(f) = \lcz(gh)$,
and since $\lexpz(f) = \lexpz(gh)$, it follows that $\lexpz(f-gh) < \lexpz(f)$.
\end{enumerate}
\end{alg}

To compute syzygy families, we first need the following lemma.

\begin{lem}\label{lem:bigenoughidealmem}
Let $M$ be a monoid such that $\Z \subset M \subset \Q$, and let $q$ be an element of the
quotient group of $M$ (i.e., the set of differences of elements of $M$).
Then for $n \gg 0$, $q + n \in M$.
\end{lem}

We now prove that the intersection of principal ideals in ${\Lambda_{}}$, both generated by elements of $\Omega_i$, must
be finitely generated by elements of $\Omega_i$.

\begin{lem}\label{lem:pip}
Given $f, g \in \kx^*$ such that $\lexpz(f), \lexpz(g) \in \Omega_i$,
there exists a finite subset of $\Omega_i$ that generates
 $\langle \lexpz(f) \rangle \cap \langle \lexpz(g) \rangle$.
 \end{lem}

\begin{proof}
By Lemma \ref{lem:bigenoughidealmem},
for each element $\sigma$ of $\Omega_i$,
 there exists a minimal
$\eta_\sigma \in \Z$ such that
$\sigma - \lexpz(f) + \eta_\sigma$, $\sigma - \lexpz(g) + \eta_\sigma \in {\Lambda_{}}$;
that is, $\sigma + n_\sigma \in \langle \lexpz(f) \rangle \cap \langle \lexpz(g) \rangle$.
Define $\Upsilon_i$ to be the finite collection  $\{ \sigma + \eta_\sigma \mid \sigma \in \Omega_i\}$.
We will show that $\Upsilon_i$ generates
 $\langle \lexpz(f) \rangle \cap \langle \lexpz(g) \rangle$.

Let $m \in \langle \lexpz(f) \rangle \cap \langle \lexpz(g) \rangle$.
By Theorem \ref{Mzunique},
${\Lambda_{}} =  \bigcup_{j=0}^\infty  \Omega_j$, and so
for some index $I$, there exist $\alpha_f, \alpha_g \in \Omega_I$
such that $m = \lexpz(f) + \alpha_f = \lexpz(g)  + \alpha_g$.
Write $\alpha_f$ as $\alpha_f' + \sum_{j=i+1}^I d_j \rho_j$ and $\alpha_g$
as $\alpha_g' + \sum_{j=i+1}^I d_j' \rho_j$ where $\alpha_f', \alpha_g' \in \Omega_i$
and $0 \le d_j, d_j' < s_j$.  By Corollary \ref{omegaclosed}, $\lexpz(f) +  \alpha_f', \lexpz(g) + \alpha_g' \in \Omega_i$.
By the uniqueness of representation promised by Theorem \ref{Mzunique},
since $m = (\lexpz(f) +  \alpha_f') + \sum_{j=i+1}^I d_j \rho_j =  (\lexpz(g) +  \alpha_g')
+ \sum_{j=i+1}^I d_j' \rho_j$, we have $d_j = d_j'$ for $i+1 \le j \le I$.
Thus $\lexpz(f) +  \alpha_f' =  \lexpz(g) +  \alpha_g'$.  So by
Theorem \ref{Mzunique},
$m':= \lexpz(f) + \alpha_f' =  \lexpz(g) + \alpha_g' =n + \sum_{j=1}^i \delta_j \rho_j$, where $n\in \N$ and $0 \le \delta_j < s_j$.
Define $\sigma = \sum_{j=1}^i \delta_j \rho_j$, and let $n_\sigma$ be the smallest
$n_\sigma \in \Z$ such that $\sigma + n_\sigma \in \langle \lexpz(f) \rangle
\cap \langle \lexpz(g) \rangle$.
Since $m' = \sigma + n \in \langle  \lexpz(f) \rangle
\cap \langle \lexpz(g) \rangle$, it follows that $n\
\ge n_\sigma$.
Thus $m' = (n-n_\sigma) + (\sigma + n_\sigma) \in \N + \Upsilon_i$,
and so
$m = m' +  \sum_{j=i+1}^I d_j \rho_j  =
(n-n_\sigma) + (\sigma + n_\sigma) +   \sum_{j=i+1}^I d_j \rho_j \in \N + \Upsilon_i + {\Lambda_{}} =
\Upsilon_i + {\Lambda_{}}$.
\end{proof}

 The following algorithm uses the lemma above to produce a syzygy family for a pair of polynomials.

\begin{alg}\label{alg:syzygy}
Let $f, g \in k[x,y]$.  This algorithm will produce
$m_1, \dots, m_\ell \in {\Lambda_{}}$ such that
$\langle \lexpz(f) \rangle \cap
\langle \lexpz(g) \rangle = \langle m_1, \dots, m_\ell \rangle$.
In addition $a_j,b_j \in \kx$ will be produced such that
$\lexpz(a_jf-b_jg) < m_j$ for each $1 \le j \le \ell$.

\begin{enumerate}
\item[(1)] Using Algorithm \ref{alg:decomposem}, write $\lexpz(f) = n + \sum_{j=1}^i d_j \rho_j$ and $\lexpz(g) = n' + \sum_{j=1}^i
d_j' \rho_j$ where $n, n'\in \N$ and $0 \le d_j, d_j' < s_j$.

\item[(2)]  Let $\sigma_1, \dots, \sigma_\ell$ be the elements of $\{\sum_{j=1}^i d_j \rho_j \mid 0 \le d_j < s_j\}$.
For each $1 \le t \le \ell$, find a minimal $\eta_t$
 such that $\sigma_t - \lexpz(f) + \eta_t, \sigma_t - \lexpz(g) + \eta_t \in {\Lambda_{}}$.
To accomplish this, begin with $\eta = 0$ and keep incrementing $\eta_t$ until
$\sigma_t - \lexpz(f) + \eta_t, \sigma_t - \lexpz(g) + \eta_t \in {\Lambda_{}}$ by
Algorithm
\ref{alg:decomposem}.

\item[(3)]  For each $t$, define $m_t = \eta_t + n_t$.  By
Lemma \ref{lem:pip}, $\{ m_1, \dots, m_\ell\}$ generates $\langle \lexpz(f) \rangle \cap
\langle \lexpz(g) \rangle$.
\end{enumerate}

\end{alg}

Below is an example of a generalized
Gr\"obner basis with respect to a valuation that is not  a Gr\"obner
basis with respect to
any monomial order.

\begin{example}
Let $k$ be a field that is not of characteristic two. Define $f_1
= y^2-x$ and $f_2 = xy$.  Then one can check that the set $B= \{
f_1,\, f_2\}$ is a Gr\"obner basis for the ideal $I = \langle f_1,
f_2 \rangle$ with respect to the valuation induced by $z = t^{1/2}
+ t^{1/4} + t^{1/8} + t^{1/16} + \cdots$ using Algorithm \ref{alg:gengbconstruct}.

We now demonstrate that $B$ is not a Gr\"obner basis with respect
to any monomial order.  Suppose, for contradiction, that $B$ is a
Gr\"obner basis with respect to some monomial order `$<$'. Note that
$x^2, y^3 \in I$ since $x^2 = yf_2 - xf_1$ and $y^3 =
y f_1 + f_2$.  We consider two cases, depending on
whether $x>y^2$ or $x<y^2$. If $x<y^2$, then lt$(f_1)=y^2$ and
lt$(f_2)=xy$.  However, $x^2 \in I$, and so if $B$ were a
Gr\"obner basis with respect to `$<$', then either $y^2 \mid x^2$ or $xy \mid x^2$, a
contradiction.  Now suppose $x>y^2$, in which case lt$(f_1)=x$ and
lt$(f_2)=xy$. However, $y^3 \in I$, and so if $B$ were a
Gr\"obner basis, then either $x \mid y^3$ or $xy \mid y^3$, a
contradiction.
\end{example}

Lastly, we note by example that some ideals do not have finite Gr\"obner bases
with respect to a given valuation.  We first prove a short lemma.

\begin{lem}\label{lem:rhoincrease}
The sequence $\rho_0, \rho_1, \rho_2, \dots$ is increasing.
\end{lem}

\begin{proof}
Since $s_j > 1$ for each index $j$, by Lemma \ref{lem:difflams},
$
\rho_i = \sum_{j=1}^{i-1}(s_j-1)\rho_j + e_{l(i)}
> \sum_{j=1}^{i-1} \rho_j + e_{l(i)} > \rho_{i-1}.
$\end{proof}

\begin{example}
Consider the ideal $\langle x,y \rangle$ of $k[x,y]$, and let $G$ be a Gr\"obner
basis with respect to the series $z\in \kq$.
For each $\rho_i$, let $p_i(x,y) \in k[x,y]$ such that
$\lexpz(p_i) = \rho_i$.
Since $G$ is a Gr\"obner basis, there exists $g_i \in G$ such that
$\lexpz(g_i) \mid \lexpz(p_i)$.
 That is, for some $h_i \in k[x,y]$, $\lexpz(g_i h_i) =
\rho_i$.  Since $G \cap k = \emptyset$, $\lexpz(g_i) > 0$,
and so $\lexpz(h_i) < \rho_i$.
Suppose, for contradiction, $\lexpz(g_i) \not= \rho_i$.
Then  $\lexpz(g_i) < \rho_i$, and so by Theorem \ref{Mzunique} and Lemma \ref{lem:rhoincrease},
 $\lexpz(g_i) = n + \sum_{j=1}^{i-1} {d_j} \rho_j$ and
 $\lexpz(h_i) = n' + \sum_{j=1}^{i-1} {d_j'} \rho_j$.
Thus, $\rho_i = \lexpz(g_ih_i) \in (1/r_{l(i-1)}) \Z$, which contradicts
Lemma \ref{cor:rhoresidue}.  Therefore, $\lexpz(g_i) = \rho_i$, and thus $G$ is infinite.
\end{example}

\end{document}